\documentclass{article}
\usepackage{amsmath}
\usepackage[left=3.09cm, right=3.09cm, bottom=2.71cm, top=2.71cm]{geometry}
\usepackage{hyperref}
\usepackage{cleveref}
\usepackage{amssymb}
\usepackage{mathtools}
\usepackage{tikz-cd}
\usepackage{amsthm}
\usepackage{graphicx}
\usepackage{graphics}
\usepackage{mathrsfs}
\usepackage{cleveref}
\usepackage{thmtools}
\usepackage{enumitem}
\usepackage{bbm}
\newlist{steps}{enumerate}{1}
\usepackage{mathtools}
\makeatletter
\newcommand{\xMapsto}[2][]{\ext@arrow 0599{\Mapstofill@}{#1}{#2}}
\def\Mapstofill@{\arrowfill@{\Mapstochar\Relbar}\Relbar\Rightarrow}
\makeatother
\usepackage{array}
\setlist[steps, 1]{label = Step \arabic*:}
\theoremstyle{plain}
\usepackage{footnote}
\makesavenoteenv{tabular}
\newtheorem*{theorem*}{Theorem}
\newtheorem{theorem}{Theorem}[subsection]
\newtheorem{definition}[theorem]{Definition}

\newtheorem{proposition}[theorem]{Proposition}
\newtheorem{remark}[theorem]{Remark}

\newtheorem*{corollary*}{Corollary}
\newtheorem{exmp}[theorem]{Example}

\newtheorem{lemma}[theorem]{Lemma}
\newtheorem*{lemma*}{Lemma}

\numberwithin{equation}{subsection}
\usepackage{rotating}
\usepackage{comment}
\usepackage[nottoc,notlot,notlof]{tocbibind} 

\usepackage[toc,page]{appendix}

\title{A Derivation of Geometric Quantization via Feynman's Path Integral on Phase Space} 
\author{Joshua Lackman\footnote{josh@pku.edu.cn} }
\date{}

\begin{document}

\maketitle
\begin{abstract}
\noindent We derive the geometric quantization program of symplectic manifolds, in the sense of both Kostant–Souriau and Weinstein, from Feynman's path integral formulation on phase space. The state space we use contains states with negative norm and polarized sections determine a Hilbert space. We discuss ambiguities in the definition of path integrals arising from the distinct Riemann sum prescriptions and its consequence on the quantization of symplectomorphisms.
\end{abstract}
\tableofcontents
\section{Introduction}
In \cite{Lackman3}, \cite{Lackman4}, the author discussed a direct relationship between Kontsevich's star product \cite{kontsevich} and the geometric quantization programs of Kostant–Souriau (\cite{bates}, \cite{woodhouse}) and Weinstein\footnote{This program emphasizes the $C^*$–algebra over the Hilbert space and uses symplectic groupoids.} (\cite{eli}, \cite{weinstein}), based on the Poisson sigma model (\cite{bon}, \cite{catt}). This resulted in a path integral definition of the quantization map, which could be put on a lattice. In this note, we will derive the same quantization map from Feynman's construction of quantum mechanics via the phase space path integral (\cite{klein}, \cite{witten}). The operators will be obtained from the expectation values of observables, as determined by Feynman's path integral. We will show that this prescription agrees with Kostant–Souriau's prequantization map on the small subspace of observables for which the latter is defined, and furthermore, that it determines sections of the prequantum line bundle over the symplectic groupoid, as in Weinstein's program.
\\\\We consider a larger space of states than usual, ie. we consider all sections of the prequantum line bundle — we find that there are states with negative norm. Polarizations are used to produce a Hilbert space (ie. states with positive norm). This gives a mathematical reason for polarizations. On the other hand, non–polarized states determine mixed states.
\\\\In addition, we discuss ambiguities in the definition of path integrals arising from the distinct Riemann sum prescriptions, ie. left-point and midpoint, and explain why this ambiguity reduces the group of symplectomorphisms which can be quantized. We suggest that this resolves a conundrum discussed in section 2.3 of \cite{brane}.
\\\\Relationships between path integrals and the Hilbert space formalism are old, and includes a relation between Feynman's path integral and Toeplitz quantization, \cite{charles}. The perspective we present is intertwined with the Poisson sigma model and symplectic groupoids. Our path integral calculations are formal, but we expect that they will be made rigorous. We will work out an explicit example in \cref{exam}. 
\section{Geometric Quantization From Feynman's Path Integral}
\subsection{Quantization of Phase Space with States of Negative Norm}
Consider a symplectic manifold $(M,\omega)$ of dimension $2n$ and let $(\mathcal{L},\nabla,\langle\cdot,\cdot\rangle)\to M$ be a prequantization. We can \textit{formally} quantize the classical theory by considering the amplitude
\begin{equation}\label{kernel}
   \langle m_1;T\,| m_0;0\rangle = \int_{\gamma(0)=m_0}^{\gamma(T)=m_1}\mathcal{D}\gamma\,e^{\frac{i}{\hbar}\int_{0}^{T}\gamma^*\nabla-\gamma^*H\,dt}
\end{equation}
for a Hamiltonian $H:M\to\mathbb{R},$ where we have implicitly chosen trivializations of $\mathcal{L}$ over $m_0,m_1\in M.$\footnote{Related is the discussion in section 2.3 of \cite{brane}, which we will comment on later.} This is \textit{formally} the amplitude for a particle, that is measured  at $m_0$ at $t=0,$ to be measured at $m_1$ at time $t=T.$ Here, the integral is over all paths connecting $m_0, m_1.$ We can compute correlation functions by inserting functions of $M$ into the integrand, eg. we can insert $f(\gamma(1/2))$ for a function $f:M\to\mathbb{C}.$
\\\\More generally, given two states $\Psi_0,\Psi_1$ (ie. sections of $\mathcal{L}$), the time evolution is determined by the phase space path integral
\begin{equation}\label{pair}
    \langle \Psi_1;T\,|\Psi_0;0\rangle=\int_{M\times M}\omega^n\boxtimes\omega^n\int_{\gamma(0)=m_0}^{\gamma(T)=m_1}\mathcal{D}\gamma\,\big\langle \Psi_1(m_1),\,e^{\frac{i}{\hbar}\int_{0}^{T}\gamma^*\nabla-\gamma^*H\,dt}\Psi_0(m_0)\big\rangle\;.
\end{equation}
Here, 
\begin{equation}
    e^{\frac{i}{\hbar}\int_{0}^{T}\gamma^*\nabla} \Psi_0(m_0)
    \end{equation}
parallel transports $\Psi(m_0)$ over $\gamma,$ and $\langle \cdot,\cdot\rangle$ in the integrand is the Hermitian inner product of $\mathcal{L}.$
\\\\Setting $H=0$ results in the inner product of the states $\Psi_0, \Psi_1$ as there is no time evolution. Explicitly,
 \begin{equation}\label{inner}
    \langle \Psi_1\,|\Psi_0\rangle=\int_{M\times M}\omega^n\boxtimes\omega^n\int_{\gamma(0)=m_0}^{\gamma(1)=m_1}\mathcal{D}\gamma\,\big\langle \Psi_1(m_1),\,e^{\frac{i}{\hbar}\int_{0}^{1}\gamma^*\nabla}\Psi_0(m_0)\big\rangle\;.
\end{equation}
\textit{Formally}, this solves the quantization problem. The problem is that this inner product space has states of negative norm — that $\langle \Psi |\Psi\rangle\in\mathbb{R}$ can be observed by noting that $\langle\cdot |\cdot\rangle$ does have conjugation symmetry, which is due to the fact that taking $i\to -i$ in the exponent is equivalent to reversing the orientation of the curve. The presence of states of negative norm is an issue that is encountered in gauge theory. Here, the existence of states of negative norm is surely related to the fact that a particle's position in phase space is not measurable, ie. $\hat{p}, \hat{q}$ don't commute. That there are states of negative norm can be checked explicitly, in the example of \cref{exam}.
\subsection{Obtaining a Hilbert Space}
One solution to the problem of negative norm states is to restrict \cref{pair} to a subspace containing only states of positive norm. This is normally (implicitly) implemented by choosing a (complex) Lagrangian polarization and taking the Hilbert space to be the subspace of polarized states.\footnote{A polarization isn't the only way of obtaining a Hilbert space, which is good because they don't always exist, \cite{go}.} This is somewhat as in gauge theory, where the Hilbert space is taken to consist of states which are constant along the leaves of a foliation determined by a group action.\footnote{In some cases, a Lagrangian polarization is also determined by a group action.} That the resulting inner product is positive definite can be seen by observing that, if $L_0, L_1$ are two Lagrangian fibers, then a path integral representation of the delta function (after possibly normalizing by an infinite constant) implies that
\begin{align}
      \int_{\gamma(0)\in L_0}^{\gamma(1)\in L_1}\mathcal{D}\gamma\,\big\langle \Psi_1(\gamma(1)),\,e^{\frac{i}{\hbar}\int_{0}^{1}\gamma^*\nabla}\Psi_0(\gamma(0))\big\rangle
    =\begin{cases}\big\langle \Psi_1(m),\Psi_0(m)\big\rangle & \textup{if } L_1=L_0\\
    0 & \textup{otherwise,}
    \end{cases}
\end{align}
where $m$ is any point in $L_0.$ This implies:
\begin{lemma}\label{lemma}(physics lemma)
On polarized sections, the inner product \cref{inner} agrees with the $L^2$–inner product of sections of the line bundle, which is given by integrating the pointwise inner product with respect to the volume form $\omega^n$ (after possibly normalizing or breaking the invariance along the polarization).    
\end{lemma}
\subsection{The Quantization Map}
Now we restrict to polarized states. To determine the expectation value of $f:M\to\mathbb{C}$ we insert $f(\gamma(1/2))$ into the integrand of \cref{inner} (it doesn't matter which interior point of the interval we evaluate $\gamma$ at, since to evaluate equal time expectation values we set $H=0$). Doing this, we see that the corresponding quantum operator $Q_f$ in the Hilbert space formalism is given by
\begin{equation}\label{hilbert}
    (Q_f\Psi)(m)=\int_{\gamma(1)=m} \mathcal{D}\gamma\, f(\gamma(1/2))\,e^{\frac{i}{\hbar}\int_0^1\gamma^*\nabla}\Psi(\gamma(0))\;.
\end{equation}
In \cite{Lackman4}, this quantization map is derived from the perspective of the Poisson sigma model.
\subsubsection{Comparison with Kostant–Souriau's Prequantization Map}
On the other hand, the Kostant–Souriau prequantization map is given by
\begin{equation}
    \frac{\hbar}{i}\nabla_{X_f}+f\;.
    \end{equation}
This map preserves polarized sections, for a polarization $\mathcal{P},$ if and only if 
\begin{equation}\label{ham}
    [X_f,\mathcal{P}]\in\mathcal{P}\;,
    \end{equation}
where $X_f$ is the Hamiltonian vector field of $f$ (this is a small space of functions, which doesn't even contain the Hamiltonian of the free particle). 
\\\\The Darboux theorem for fibrations \cite{arnold} states that any Lagrangian polarization is locally given by constant $q$-coordinates in a canonical coordinate system $(p_1,\ldots,p_n,q^1,\ldots,q^n).$ In such a coordinate system, $f$ satisfies \cref{ham} if and only if it is affine in each $p_i.$ On such functions, that the Kostant–Souriau prequantization map gives the correct operator (ie. \cref{hilbert}) can be seen from a local computation, as done in \cite{klein}, section 2.1.2. Essentially, since $\frac{\partial}{\partial q^i}$ is the Hamiltonian vector field of $p_i,$ it is due to the equation
\begin{equation}
    p_i(t)\,e^{\frac{i}{\hbar}\int_0^1 p_1(t')\frac{dq^1(t')}{dt'}+\cdots p_n(t')\frac{dq^n(t')}{dt'}\,dt'}=\frac{\hbar}{i}\frac{\delta}{\delta q^i(t)}e^{\frac{i}{\hbar}\int_0^1 p_1(t')\frac{dq^1(t')}{dt'}+\cdots p_n(t')\frac{dq^n(t')}{dt'}\,dt'}
\end{equation}
and integration by parts.
\subsubsection{Quantized Observables as States of the Fundamental Groupoid}
We can think of $Q_f$ as an integral kernel, in the following way: let $\Pi_1(M)$ denote the manifold which consists of homotopy classes of maps $\gamma:[0,1]\to M$ relative to the endpoints, and denote such an equivalence class by $[\gamma].$ There are two maps 
\begin{equation}
    s,t:\Pi_1(M)\to M
    \end{equation}
given by evaluation at $0,1,$ respectively. Thus, we can form the line bundle
\begin{equation}\label{line}
    t^*\mathcal{L}\otimes s^*\mathcal{L}^*\to \Pi_1(M)\;,
\end{equation}
where $\mathcal{L}$ is the prequantum line bundle of $(M,\omega).$ A point in this line bundle over $[\gamma]$ defines a linear map between the fibers of $\mathcal{L}\to M$ over $\gamma(0),\, \gamma(1).$ If $M$ is simply connected then this space is just $M\times M$ and the points in this line bundle are exactly identified with such linear maps. Otherwise 
\begin{equation}
    \Pi_1(M)=\widetilde{M}\times \widetilde{M}/\pi_1(M)\;,
    \end{equation}
where $\widetilde{M}\to M$ is the universal cover and where $\pi_1(M)$ acts according to the diagonal action. In this case, distinct points in the line bundle can define the same linear map. Now, the point is that $Q_f$ defines a section of \cref{line}:\footnote{We are abusing notation in identifying the operator $Q_f$ with this section of the line bundle, but it is justified by \cref{just}.}
\begin{equation}
    Q_f([\gamma'])=\int_{\gamma\in[\gamma']}\mathcal{D}\gamma\,f(\gamma(1/2))\,e^{\frac{i}{\hbar}\int_0^1\gamma^*\nabla}\;.
\end{equation}
We can then rewrite \cref{hilbert} as 
\begin{equation}\label{just}
    (Q_f\Psi)(m)=\int_{[\gamma]\in t^{-1}(m)} [s^*\omega^n] \;Q_f([\gamma])\Psi(s[\gamma]))\;.
\end{equation}
Note that, in the above formula 
\begin{equation}
    s:t^{-1}(m)\to M
    \end{equation}
is identified with the universal cover with basepoint $m,$ and $s^*\omega^n$ is the induced measure. On a symplectic manifold\footnote{This can be generalized to a Poisson manifold.} this gives a definition of the quantization map in Weinstein's program. Functions on the phase space which are square integrable should result in genuine sections of the line bundle over $\Pi_1(M),$ whereas other functions may result in distribution–valued sections. Note that, the line bundle \cref{line} is the prequantum line bundle with respect to $t^*\omega-s^*\omega,$ making it a symplectic groupoid.
\\\\Higher goupoids provide the natural mathematical setting for making these constructions. For a symplectic manifold the relevant (higher) groupoids are the fundamental groupoid, which we have denoted $\Pi_1(M),$ and the fundamental 2-groupoid, denoted $\Pi_2(M),$ which is the truncation in degree $2$ of the singular simplicial space of $M.$ The 1-arrows in $\Pi_2(M)$ are paths in $M$ — which is the domain of integration of the path integral. The 2-arrows are homotopies of paths, relative to the endpoints\footnote{For symplectic groupoids, eg. $\Pi_1(M),$ see \cite{ruif}. For the relevant Lie 2-groupoid, eg. $\Pi_2(M),$ see \cite{zhuc}.}—these are things that we can integrate $\omega$ over.
\\\\In the more general case of Poisson manifolds discussed in \cite{Lackman3}, the quantizations of observables which are supported on distinct isomorphism classes of points in phase space commute. On (connected) symplectic manifolds, all points in phase space are isomorphic (in the sense that there is an arrow connecting any two of them). Only in the case of the zero Poisson structure (ie. when $\hbar=0$) do all points represent distinct classes.
\\\\The operator composition $Q_fQ_g$ corresponds to the (twisted) convolution $Q_f\ast Q_g$ of $\Pi_1(M),$\footnote{The twisted convolution is the convolution of functions valued in a line bundle with a multiplication. See \cite{eli} for the definition.} ie.
\begin{equation}
Q_fQ_g\Psi=(Q_f\ast Q_g)\Psi\;.
\end{equation}
This is \textit{very} closely related to the Poisson sigma model formulation of Kontsevich's star product \cite{kontsevich}, which is \textit{implicitly} argued in \cite{bon} to be a twisted convolution product in $\Pi_2(M)$ (shown in \cite{Lackman3}, \cite{Lackman4}). The star product on symplectic manifolds can be written as
\begin{equation}\label{moy}
    (f\star g)(m)=\int_{\gamma(\infty)=m} \mathcal{D}\gamma\, f(\gamma(0))g(\gamma(1))\,e^{\frac{i}{\hbar}\int_{S^1}\nabla^*\gamma}\;,
\end{equation}
where the integral is over all contractible loops $\gamma:S^1\to M$ and $0,1,\infty$ are three cyclically ordered marked points,  \cite{brane}, \cite{grady}.
\begin{remark} $Q_f,$ as well as the product \cref{moy}, make sense as long as the pullback of $\omega$ to $\widetilde{M}$ is prequantizable. This more general approach allows us to quantize the torus with translation invariant symplectic form $\omega/\hbar\,,$ for any $\hbar,$ as done in \cite{weinstein}.
\end{remark}
\subsection{The Abstract $C^*$–algebra and Mixed States}
Importantly, the quantization map $f\mapsto Q_f$ generates an abstract $C^*$–algebra, which is a $C^*$-subalgebra of the twisted convolution algebra of $\Pi_1(M),$ \cite{eli}. A reference for this section is \cite{williams}.
\\\\A state of a unital $C^*$–algebra $\mathcal{A}$ is a normalized, positive linear functional, ie.
\begin{definition}
A state of a unital $C^*$–algebra $\mathcal{A}$ is a linear functional $\rho:\mathcal{A}\to\mathbb{C}$ such that $\rho(1)=1$ and $0\le \rho(A^*A)\in\mathbb{R}$ for all $A\in\mathcal{A}.$ States that can't be written as a convex linear combination of distinct states are called pure states, the other states are called mixed states.\footnote{Pure states are states of maximal information, mixed states are statistical ensembles.}
\end{definition}
Let us observe that the definition of $Q_f$ makes sense independent of any polarization, and thus it acts on the entire prequantum Hilbert space. For this reason, we can compute expectation values of any section of the prequantum line bundle, with respect to the $L^2$–inner product of sections. Therefore, states which aren't polarized naturally acquire the interpretation of mixed states, while polarized states are pure states. 
\\\\To say a bit more, there is a simple correspondence between states of a $C^*$–algebra and its representations as operators on a Hilbert space, for which pure states correspond to irreducible representations. A vector in a representation is called a vector state, and its expectation values determine a state. This gives the correspondence in one direction. The correspondence in the other direction is given by the GNS construction:
\\\\\textit{Let $\rho$ be a state of $\mathcal{A}$ and consider the quotient $\mathcal{A}/\sim\,,$ where we quotient out by the left ideal consisting of those $A\in\mathcal{A}$ such that $\rho(A^*A)=0.$ There is a natural inner product on $\mathcal{A}/\sim\,,$ given by 
\begin{equation}
    \langle A,B\rangle:=\rho(A^*B)\;.
    \end{equation}
The completion is a Hilbert space and it forms a representation of $\mathcal{A},$ where $\langle 1,A1\rangle=\rho(A),$ so that the state determined by $1\in\mathcal{A}/\sim$ is equal to $\rho.$}
\\\\Suppose that the $C^*$–subalgebra generated by $\{A_1,\ldots, A_n\}\subset \mathcal{A}$ is commutative. Then the pullback of a pure state to this $C^*$–subalgebra is a mixed state that determines a probability measure on a topological space, via Gelfand duality. The standard interpretation is that this measure gives the probability of measuring $A_1,\ldots, A_n$ to have given (simultaneous) values, which correspond to points in the topological space — this $C^*$–subalgebra being commutative is equivalent to $A_1,A_1^*,\ldots,A_n,A_n^*$ pairwise commuting. Of course, if an operator commutes with its adjoint it is measurable since its real and imaginary parts are self-adjoint and commute.
\section{Computing the Path Integrals: Quantum Mechanics on the Cotangent Bundle}\label{exam}
\subsection{Computing the Inner Product}
Now we will discuss \cref{kernel} on $\textup{T}^*\mathbb{R}$ in the case that $H=0,$ as this is needed to compute the inner product and quantization map. To do this we should first prequantize $\textup{T}^*\mathbb{R}\cong \mathbb{R}^2.$ There is a standard construction of such a prequantization: since $\textup{T}^*\mathbb{R}$ is simply connected, any two paths relative to the endpoints are homotopic. Let $D$ be a disk with marked north and south pole. Consider a map
\begin{equation}
    X:D\to \textup{T}^*\mathbb{R}\,,
\end{equation}
and denote the restriction of $X$ to the left and right hemispheres by 
\begin{equation}
    \gamma_1,\,\gamma_2:[0,1]\to \textup{T}^*\mathbb{R}\;.
\end{equation}
We identify 
\begin{equation}
    (\gamma_1,\lambda\sim(\gamma_2,\lambda\,e^{\frac{i}{\hbar}\int_D X^*\omega})\in )\in C([0,1],\textup{T}^*\mathbb{R})\times\mathbb{C}\;.
\end{equation}
Doing this for all $X$ determines a line bundle with Hermitian connection over $\textup{T}^*\mathbb{R}\times \textup{T}^*\mathbb{R}.$ We can pull back this line bundle to $\{(0,0)\}\times \textup{T}^*\mathbb{R},$ and this determines a line bundle with connection over $\textup{T}^*\mathbb{R}.$ That is, we are only considering paths starting at $(0,0).$ Next, we will describe this line bundle more explicitly.
\\\\Choose a complex structure turning $\textup{T}^*\mathbb{R}\cong \mathbb{R}^2$ into a K\"{a}hler manifold. This determines a natural trivialization of the line bundle over $\mathbb{R}^2,$ since between any two points is a unique geodesic, and therefore for any $(p,q)\in \mathbb{R}^2$ there is a unique geodesic connecting $(0,0)$ to $(p,q).$ 
\\\\The resulting line bundle with Hermitian connection is
\begin{equation}
    \mathcal{L}=\mathbb{R}^2\times\mathbb{C}\to\mathbb{R}^2\;,\;\nabla=\frac{pdq-qdp}{2\hbar}\;.
    \end{equation}
We then have that (up to a finite normalization constant)
\begin{equation}\label{ans}
   \langle p_1,q_1\,| p_0,q_0\rangle = \int_{\gamma(0)=(p_0,q_0)}^{\gamma(1)=(p_1,q_1)}\mathcal{D}\gamma\,e^{\frac{i}{2\hbar}\int_{0}^{1}\gamma^*(pdq-qdp)}=e^{\frac{i}{2\hbar}(p_0q_1-q_0p_1)}\;.\footnote{See appendix B of \cite{Lackman3} for a discussion about the computation.}
\end{equation}
The quantity on the right side is equal to the parallel transport of $1$ from $(p_0,q_0)$ to $(p_1,q_1)$ over the corresponding \textit{geodesic} (or really, it's pairing with $1$ over $(p_1,q_1)).$ Before continuing this example, we should address the fact that \cref{ans} seems to depend on the choice of a metric, which doesn't appear in \cref{inner}.
\subsection{A Conundrum of the Propagator}\label{conun}
A Riemannian structure doesn't appear in \cref{inner}, so the appearance of a geodesic on the right side of \cref{ans} may be surprising, since we would get a different answer if we changed the complex structure. This is related to a conundrum discussed in section 2.3 of \cite{brane}, where it is suggested that due to the existence of too many symmetries, if the path integral is defined then it must be the case that 
\begin{equation}
    \langle p_1,q_1\,| p_0,q_0\rangle=\delta(p_1-p_0,q_1-q_0)\;,
    \end{equation}
which contradicts \cref{ans}. However, as we will see when we continue our example in \cref{continue}, \cref{ans} does produce the correct quantum mechanics in the context of \cref{lemma} (it also produces the correct operators). The resolution of this conundrum is that the symmetries of the path integral \cref{inner} aren't all symplectomorphisms. This is for the same reason that, when discretizing the path integral of a particle in a magnetic potential $A$ (\cite{klein}), the left-point and midpoint Riemann sums corresponding to 
\begin{equation}
    \int_0^1 A(q)\,\frac{dq}{dt}dt
\end{equation}
give different results for the path integral, in the limit as the lattice spacing goes to zero (which is related to the operator ordering problem). That is, this path integral isn't \textit{a priori} well-defined because there's a choice involved in which Riemann sum prescription is used. Essentially, different Riemann sum prescriptions preserve different symmetries, and since the limits don't agree, we don't retrieve the full group of symplectomorphisms. This is also the reason why there is no canonical quantization map for a symplectic manifold. However, on K\"{a}hler manifolds there is a canonical prescription for the Riemann sums — this is consistent with Fedesov's construction of the star product \cite{fedosov}, which is canonical on K\"{a}hler manifolds, as well as Berezin–Toeplitz quantization \cite{bord}.
\\\\This phenomenon of different Riemann sum prescriptions resulting in different limits is also the reason that the It\^{o} and Stratonovich integrals on Wiener space are distinct, even though classically they are the same, as the former is computed using a left-point Riemann sum and the latter is computed using the midpoint Riemann sum. It's a consequence of the fact that the paths of the path integral aren't smooth, and so higher order terms appearing in a Riemann sum give nontrivial contributions to the resulting integral. We will discuss this more precisely next. For the basics of Wiener space (Brownian motion) and stochastic integrals, see \cite{rick}, \cite{bernt}.
\subsection{Formalizing Riemann Sums Under the Path Integral}
For a smooth function 
\begin{equation}
    f:\mathbb{R}\to\mathbb{R}\;,
\end{equation}
the left-point and midpoint evaluations
\begin{equation}
    f(x_i)\,, \;f\Big(\frac{x_i+x_{i+1}}{2}\Big)
\end{equation}
differ by a term of order $(\Delta x_i)^2,$ whereas only terms of order $\Delta x_i,$ in a traditional Riemann sum, contribute to the integral. If we integrate $f\,dx$ over a \textit{smooth} path 
\begin{equation}
    x:[0,1]\to \mathbb{R}\;,
    \end{equation}
it doesn't matter which prescription we use for the Riemann sums, since for such paths $\Delta x^2\sim\Delta t^2.$\footnote{Here, $\Delta x^2$ means $(\Delta x)^2.$} However, the paths of the path integral are \textit{rough}, and they are such that 
\begin{equation}
    \Delta x^2\sim \Delta t\;.
    \end{equation}
To make a precise statement, the paths in Wiener space (ie. the Euclidean path integral) have H\"{o}lder continuity with exponent less than $1/2,$ almost surely (by L\'{e}vy's modulus of continuity, \cite{peter}). Therefore, the higher order terms do contribute to the integral and the different prescriptions for the Riemann sums have different limits. This is essentially the cause of what physicists call the ``operator ordering problem". To be mathematically precise:
\begin{proposition}\label{wien}\footnote{The Riemann sum on the left defines the  Stratonovich integral and the Riemann sum on the right defines the It\^{o} integral, \cite{bernt}. The former satisfies the fundamental theorem of calculus.}
With respect to the Wiener measure on continuous paths $x\in C([0,1],\mathbb{R})\,,$ and for a smooth function $f:\mathbb{R}\to\mathbb{R}\,,$ 
\begin{align}
    \nonumber &\sum_{i=0}^{n-1}f\bigg(\frac{x(t_i)+x(t_{i+1})}{2}\bigg)\,(x(t_{i+1})-x(t_i))\,-\,\sum_{i=0}^{n-1}f(x(t_i))\,(x(t_{i+1})-x(t_i)) 
    \xrightarrow[]{\Delta t_i\to 0}\,\frac{1}{2}\int_0^1 \frac{df}{dx}(x(t))\,dt
    \end{align}
in $L^2(C[0,1[).$
\end{proposition}
Given this, it can be said that the things that are Riemann integrated over paths under the quantum mechanical path integral (in the Lagrangian formulation) aren't locally of the form $f\,dx,$ they are of the form 
\begin{equation}
    f\,dx+g\,dx^2\;.
\end{equation} 
This can be interpreted as a symmetric tensor or as a pointwise polyomial on the tangent spaces. The term $g\,dx^2$ is irrelevant classically since it doesn't contribute when integrated over a smooth path. This is made more precise in \cite{Lackman2}, section 6.1.
\\\\To clarify exactly how the Riemannian metric enters into \cref{ans} and in what sense the terms on the left of \cref{wien} are Riemann sums,\footnote{The terms look like Riemann sums if the integral were over $x,$ however, it is $t$ that is being integrated over.} it helps to reformulate the Riemann integral in a more mathematically natural way, which doesn't depend on coordinates and behaves well with respect to pullbacks — this is important because path integrals involve pullbacks over paths. To motivate this we state the following result.
\begin{proposition}\label{int}(\cite{Lackman2}, appendix \ref{riemann})
Let $f:[0,1]\to\mathbb{R}$ be a smooth function and let $F:[0,1]\times[0,1]\to\mathbb{R}$ be a smooth function which vanishes on the diagonal, such that
\begin{equation}\label{condi}
    \partial_y F(x,y)\vert_{y=x}=f(x)\;.
\end{equation}
Then 
\begin{equation}\label{rs}
    \sum_{i=0}^{n-1}F(x_i,x_{i+1})\xrightarrow[]{\Delta x_i\to 0}\int_0^1 f\,dx\;.
\end{equation}
\end{proposition}
As an example, $F(x,y)=f(x)(y-x),\; f(y)(y-x)$ satisfy \cref{condi}, and using these on the left hand side of \cref{rs} gives the left-point and right-point Riemmann sums, respectively. Stated more geometrically, we want the exterior derivative of $F$ in the second factor, evaluated on the diagonal, to agree with $f\,dx.$ The left side of \cref{rs} can be considered to be a generalized Riemann sum, on the manifold $[0,1]$ — indeed, \cref{int} generalizes to manifolds $M$ and gives a notion of Riemann sums of differential forms on manifolds, as described in appendix \ref{riemann}.
\begin{exmp}
The phase space path integral has an integrand which formally contains the term
\begin{equation}\label{cont}
    e^{\frac{i}{\hbar}\int_0^1 p(t)\frac{dq}{dt}\,dt}\;.
\end{equation}
Evaluating the phase space path integral involves first partitioning $[0,1]$ with points $0=t_0<t_1<\ldots<t_n=1$ and computing a discretized approximation, whose integrand contains terms such as 
\begin{equation}\label{dis}
    e^{\frac{i}{\hbar}\sum_{k=0}^{n-1}p_k (q_{k+1}-q_k)}\;,
\end{equation}
the exponent of which is taken to mean
\begin{equation}\label{nr}
    \sum_{i=0}^{n-1}p(t_k)(q(t_{k+1})-q(t_k))\;.
\end{equation}
\Cref{nr} isn't a Riemann sum of $t\mapsto p(t)\frac{dq}{dt}$ in the traditional sense. However, it is a Riemann sum of $p(t)\frac{dq}{dt}$ in the more general sense of \cref{int}. This is because the function
\begin{equation}
    [0,1]^2\to\mathbb{R}\,,\;\;(t,t')\mapsto p(t)(q(t')-q(t))\,,
\end{equation}
appearing in \cref{nr}, satisfies the conditions of \cref{int}.
\end{exmp}
$\,$\\Due to \cref{wien}, in the context of Wiener space $L^2(C([0,1]))$ we need to specify the second derivative of $F$ in the second factor as well. That is, to define 
\begin{equation}
    \big(t\mapsto x(t)\big)\mapsto \int_0^1 f(x(t))\,dx(t)+g(x(t))\,dx^2(t)
\end{equation}
as a function in $L^2(C[0,1]),$ we need to choose $F$ so that
\begin{equation}\label{cond}
    \partial_y F(x,y)\vert_{y=x}=f(x)\,,\; \frac{1}{2}\partial_y^2 F(x,y)\vert_{y=x}=g(x)\;.
\end{equation}
The result is well-defined. That is, if $F_1, F_2$ both vanish on the diagonal and satisfy \cref{cond}, then in $L^2(C([0,1]))$
\begin{equation}
      \lim\limits_{\Delta t_i\to 0}\sum_{i=0}^{n-1}F_1(x(t_i),x(t_{i+1}))=\lim\limits_{\Delta t_i\to 0}\sum_{i=0}^{n-1}F_2(x(t_i),x(t_{i+1}))\;.
\end{equation}
If we let $g=0$ then we get the It\^{o} integral.
\\\\
This explains the appearance of the metric in \cref{ans}. In order to integrate $\nabla=(pdq-qdp)/2$ in the path integral, we need to choose $F$ so that its exterior derivative in the second factor, evaluated on the diagonal, is equal to $(pdq-qdp)/2.$ A Riemannian metric determines such an $F$: we simply define $F(p_0,q_0,p_1,q_1)$ to be the integral of $(pdq-qdp)/2$ over the convex hull of $(p_0,q_0)$ and $(p_1,q_1),$ ie. the geodesic connecting them. The result is
\begin{equation}
    F(p_0,q_0,p_1,q_1)=\frac{p_0q_1-q_0p_1}{2}\;.
\end{equation}
Since the higher order terms in a Riemann sum contribute when integrated over rough paths, the path integral itself ends up depending on the metric.
\begin{exmp}
Suppose $dF=f\,dx$ on $[0,1].$ Then $(x,y)\mapsto F(y)-F(x)$ satisfies the conditions of \cref{int}, and therefore 
\begin{equation}
    F(1)-F(0)=\sum_{i=0}^{n-1}F(x_{i+1})-F(x_i)\xrightarrow[]{\Delta x_i\to 0}\int_0^1 f\,dx\;.
\end{equation}
\end{exmp}
\subsection{Continuing with the Example}\label{continue}
In order to compute \cref{hilbert}, we can rewrite the path integral as 
\begin{equation}\label{hilbert2}
    (Q_f\Psi)(m)=\int_{M\times M}\omega^n\boxtimes\omega^n\,f(m_1)\, \bigg[\int_{0,1/2,1}\mathcal{D}\gamma\, e^{\frac{i}{\hbar}\int_0^1\gamma^*\nabla}\bigg]\Psi(m_0)\;,
\end{equation}
where the integral on the right is over the set of paths with $\gamma(0)=m_0,\,\gamma(1/2)=m_1, \gamma(1)=m,$ and the first factor of $\omega^n\boxtimes\omega^n$ integrates over $m_0$ and the second integrates over $m_1.$ Absorbed into $\mathcal{D}\gamma$ on the right is a normalization factor, which is the reciprocal of
\begin{equation}
    \int_{0,1/2,1} \mathcal{D}\gamma\, e^{\frac{i}{\hbar}\int_0^1\gamma^*\nabla}\;,
\end{equation}
where here the integral is over all maps $\gamma:S^1\to M$ such that $\gamma(0)=m_0,\,\gamma(1)=m_1,\,\gamma(\infty)=m,$ where $0,1,\infty$ are three cyclically ordered marked points. Using \cref{ans}, on $\textup{T}^*\mathbb{R}$ the result is
\begin{align}\label{defi2}
    Q_f\Psi(u)=\frac{1}{(2\pi\hbar)^2}\int_{\textup{T}^*\mathbb{R}\times \textup{T}^*\mathbb{R}}\omega\boxtimes\omega\,f(v)\Psi(z)e^{\frac{i}{\hbar}\Omega(u,v,z)}\;,
\end{align}
where $\Omega(u,v,z)$ is the area of the convex hull of $u,v,z$ (ie. the corresponding triangle). This is because 
\begin{equation}
    \int_{S^1}\gamma^*\nabla=\int_{D}X^*\omega\;,
\end{equation}
where $X:D\to \textup{T}^*\mathbb{R}$ is any disk whose restriction to the boundary agrees with $\gamma.$
\\\\One can explicitly check that on states polarized with respect to (complex) linear polarizations, \cref{inner}, which is equal to 
 \begin{equation}\label{inner2}
    \langle \Psi_1\,|\Psi_0\rangle=\frac{1}{(2\pi\hbar)^2}\int_{\textup{T}^*\mathbb{R}\times \textup{T}^*\mathbb{R}}dp_0 dq_0 dp_1 dq_1\,\overline{\Psi}_1(p_1,q_1) \Psi_0(p_0,q_0)\,e^{\frac{i}{2\hbar}(p_0q_1-q_0p_1)}\;,
\end{equation}
agrees with the $L^2$–inner product of sections. Here, we are identifying $\Psi_0,\Psi_1$ with $\mathbb{C}$–valued functions (this works out nicer for complex–polarized states, since they are actually integrable with respect to $\omega$ and so we don't need to normalize). Explicitly, (complex) linearly–polarized states are of the form
\begin{equation}
\Psi(p,q)=e^{\frac{i}{2\hbar}(ap+bq)(cp+dq)}\psi(ap+bq)
\end{equation} 
for $ad-bc=1.$
\\\\Furthermore, for any state $\Psi$ polarized with respect to a (complex) linear polarization, $Q_f\Psi$ is still a polarized state. In fact, on polarized states the formula \cref{defi2} simplifies to
\begin{equation}\label{simp}
Q_f\Psi(u)=\frac{1}{2\pi\hbar}\int_{\textup{T}^*\mathbb{R}}f(u')\Psi(u'-u)e^{\frac{i}{\hbar}\Omega(0,u,u')}\,\omega\;.
\end{equation}
This formula agrees with the (non–perturbative) Weyl quantization (\cite{folland}, page 79) on states polarized with respect to the projection map, and on (complex) linearly–polarized states it agrees with the Kostant–Souriau prequantization map on the small subspace of functions for which the latter is defined. That \cref{defi2} is defined on the entire prequantum Hilbert space gives an explanation for why the BKS pairing (\cite{bates}) defines a unitary equivalence between states polarized with respect to (complex) linear polarizations.\footnote{For other polarizations, the Blattner-Kostant-Sternberg (BKS) pairing need not define a unitary map.} In particular, for any (complex) linearly–polarized states, the quantization determined by $Q_f$ is naturally equivalent to the Weyl quantization. This is discussed more in \cite{Lackman3}, \cite{Lackman4}.
\\\\In \cref{lemma}, we suggested that for any polarized states, $\langle \Psi_1\,|\Psi_0\rangle$ agrees with the $L^2$–inner product. However, this wasn't a precise statement. That in \cref{inner2} we have agreement only for (complex) linearly–polarized states can be attributed to the higher order contributions of the Riemann sums, which reduce the symmmetry group of \cref{inner} from symplectomorphisms to automorphisms of the K\"{a}hler structure.
\begin{appendices}\label{appen}
\section{Riemann Sums on Manifolds}\label{riemann}
Traditionally, Riemann sums are only defined on coordinate space. Here, we will explain (generalized) Riemann sums on manifolds, whose construction is motivated by path integrals and assigns a Riemann sum to a differential form. The simplest way to explain it uses the van Est map (originally defined for groupoids in \cite{weinstein1}), which is defined with respect to a Lie groupoid. We will describe it in the simplest case, where the relevant groupoid is the pair groupoid of a manifold $M$ (this groupoid has a unique arrow between any two points in $M$). This is discussed more in \cite{Lackman1}, \cite{Lackman2}.
\begin{definition}\label{des}\footnote{This formulation is presented in \cite{Lackman2}, section 3. It is equivalent to, but a bit different from, the original formulation in \cite{weinstein1}.}
Let $M$ be a manifold. Let $\Omega:M^{n+1}\to\mathbb{C}$ be a smooth function which vanishes on the diagonal and which is invariant under even permutations. The (graded) van Est map applied to $\Omega$ gives an $n$-form $VE_0(\Omega)$ on $M.$ Letting $X_1,\ldots X_n$ be vectors at a point $m\in M,$ it is defined by
\begin{equation}
    VE_0(\Omega)(X_1,\ldots,X_n)=n!\,X_n\cdots X_1\,\Omega(m,\cdot,\ldots,\cdot)\;,
\end{equation}
where $X_i$ differentiates $\Omega(m,\cdot,\ldots,\cdot)$ in the $ith$ component (where $m$ occupies the $0th$ component).
\end{definition}
\begin{exmp}
Let $M=\mathbb{R}^2$ and let $\Omega=\Omega(x_0,y_0,x_1,y_1,x_2,y_2)$ be a function on $\mathbb{R}^2\times\mathbb{R}^2\times\mathbb{R}^2$ which satisfies the conditions of the previous definition. Then
\begin{equation}
    VE_0(\Omega)(\partial_x\vert_{(x_0,y_0)},\partial_y\vert_{(x_0,y_0)})=2\,\partial_{y_2}\partial_{x_1}\Omega(x_0,y_0,\cdot,\cdot,\cdot,\cdot)\vert_{(x _1,y_1)=(x_2,y_2)=(x_0,y_0)}\;.
\end{equation}
\end{exmp}
Now, let $M$ be an oriented, $n$-dimensional manifold. Given a triangulation $\Delta_M$ of $M,$ up to even permutation there is a canonical ordering of the $(n+1)$ vertices of each $n$–dimensional simplex. Therefore, we can evaluate $\Omega:M^{n+1}\to\mathbb{C}$ on any $n$-dimensional face $\Delta\in\Delta_M$ by choosing such an ordering of the vertices and plugging them into $\Omega.$ This is well-defined because we are assuming that $\Omega$ is invariant under even permutations. We denote its value by $\Omega(\Delta).$
\begin{definition}(\cite{Lackman1}, \cite{Lackman2})
Let $M$ be an oriented $n$-dimensional manifold, let $\omega$ be an $n$-form on $M$ and let $VE_0(\Omega)=\omega.$ Then given a triangulation $\Delta_M$ of $M,$ the (generalized) Riemann sum of $\omega$ is defined to be
\begin{equation}
    \sum_{\Delta\in \Delta_M}\Omega(\Delta)\;,
\end{equation}
where the sum is over all $n$-dimensional simplices.
\end{definition}
\begin{theorem}\label{inte}(\cite{Lackman1}, \cite{Lackman2})
Suppose that $\Omega$ satisfies the conditions of \cref{des} and that $VE_0(\Omega)=\omega.$ Then
\begin{equation}
\sum_{\Delta \in\Delta_M }\Omega(\Delta)\xrightarrow[]{\Delta\to 0}\int_M\omega\;,
\end{equation}  
where the limit is taken over barycentric subdivisions of any triangulation $\Delta_M.$
\end{theorem}
The next proposition shows that this notion of Riemann sum is well–behaved with respect to pullbacks, which is important because pullbacks appear in path integrals:
\begin{proposition}
Let $\Delta_M$ be a triangulation of an $n$-dimensional manifold $M,$ let $f:M\to N$ be a smooth function and let $\omega$ be an $n$-form on $N$ with $VE_0(\Omega)=\omega.$ Then 
\begin{equation}
  \sum_{\Delta\in \Delta_M}f^*\Omega(\Delta)  
\end{equation}
is a (generalized) Riemann sum of $f^*\omega,$ where $f^*\Omega$ is the pullback via the map $M^{n+1}\to N^{n+1}$ induced by $f.$
\end{proposition}
This integral comes with a generalization of the fundamental theorem of calculus, and we will state one part of it. In the following, $M$ is a compact manifold with boundary, $\Omega_M$ is a function satisfying the conditions of \cref{des} on a neighborhood of the diagonal in $M^{n+1},$ and $\Omega_{\partial M}$ is a function satisfying the conditions on a neighborhood of the diagonal in $(\partial M)^n.$ Furthermore, $\Delta_{\partial M}$ is the induced triangulation on the boundary. We will state the meaning of the pair being closed after stating the result:
\begin{theorem}\label{fund}(\cite {Lackman2}, section 5)
If $(\Omega_M, \Omega_{\partial M})$ is closed as a pair of functions on $\textup{Pair}^{(n)}_{\textup{loc}}\,(M,\partial M),$ then 
 \begin{equation}\label{sing}
\sum_{\Delta\in\Delta_M}\Omega_M(\Delta)\;-\sum_{\Delta\in\Delta_{\partial M}}\Omega_{\partial M}(\Delta)
=\int_M VE_0(\Omega_M)-\int_{\partial M}VE_0(\Omega_{\partial M})\;.\footnote{This is an equality of cap products with the fundamental class of $M$: on the left is the cap product on simplicial cohomology and on the right is the cap product on de Rham cohomology.}
\end{equation}   
\end{theorem}
The meaning of the pair $(\Omega_M,\Omega_{\partial M})$\footnote{$\Omega_M$ is defined on a neighborhood of the diagonal in $M^{n+1}$ and $\Omega_{\partial M}$ is defined on a neighborhood of the diagonal in $(\partial M)^{n}.$} being closed is analogous to a pair being closed in (relative) singular cohomology: there are $(n+2)$ projection maps 
\begin{equation}
   \pi_0,\pi_1,\ldots,\pi_{n+1}: M^{n+2}\to M^{n+1}\;,
\end{equation}
where $\pi_i$ drops the $ith$ component. Similarly, there are $(n+1)$ projections maps $(\partial M)^{n+1}\to (\partial M)^{n}.$ Let $i:(\partial M)^{n+1}\to M^{n+1}$ be the inclusion. That $(\Omega_M, \Omega_{\partial M})$ is closed means that in a neighborhood of the diagonals,\footnote{These maps define a differential on functions on the nerve of $\textup{Pair}_{\textup{loc}}\,(M,\partial M).$ There is a cohomology theory using the van Est map, for functions on differentiable stacks which are valued in any abelian Lie group, \cite{Crainic}, \cite{Lackman0}, \cite{Lackman}.}
\begin{align}
&\sum_{i=0}^{n+1} (-1)^i\pi_i^*\Omega_M=0\;,
\\& i^*\Omega_M-\sum_{i=0}^n (-1)^i\pi_i^*\Omega_{\partial M}\;.
\end{align}
\begin{proof}(\cref{fund})
The proof generalizes the traditional proof of the fundamental theorem of calculus. The pair $(\Omega_M,\Omega_{\partial M})$ defines a class in the (relative) singular cohomology of $X,$ so the left side of \cref{sing} is independent of the triangulation. Taking a limit over triangulations and using \cref{inte} gives the result.
\end{proof}
To explain why it generalizes the fundamental theorem of calculus:
\begin{exmp}
Let $f:[a,b]\to\mathbb{R},$ let $F(x,y)=f(y)-f(x)$ on $[a,b]^2,$ and let $\Delta_{[a,b]}$ be the trivial triangulation of $[a,b],$ ie. the one with only one top dimensional face. Then the pair 
\begin{equation}
   (F,0)
\end{equation}
is closed as a pair on $\textup{Pair}^{(1)}_{\textup{loc}}\,\big([a,b],\{a,b\}\big).$ That is,
\begin{equation}
 F(y,z)-F(x,z)+F(x,y)=0\;,
    \end{equation}
and since we can choose a neighborhood of the diagonal of $\{a,b\}^2$ so that it only contains the points $(a,a), (b,b),$ the pullback of $F$ via the inclusion vanishes as well. Therefore, since $VE(F)=df\,,$ \cref{fund} implies that
\begin{equation}
    f(b)-f(a)=\int_a^b df\;.
\end{equation}
\end{exmp}
\end{appendices}

 \end{document}